\documentclass[12pt]{amsart}
\usepackage{amscd,amssymb,times,fullpage}
\newtheorem{theorem}{Theorem}

\newtheorem{proposition}[theorem]{Proposition}
\newtheorem{corollary}[theorem]{Corollary}
\theoremstyle{definition}
\newtheorem{definition}[theorem]{Definition}

\newtheorem{example}[theorem]{Example}

\theoremstyle{remark}
\newtheorem{remark}[theorem]{Remark}

\newcommand{\R}{\mathbb{ R}}

\newcommand{\HHH}{\mathcal{ H}}
\newcommand{\GG}{\mathcal{ G}}

\newcommand{\FF}{\mathcal{ F}}

\def\ra{{\rightarrow}}

\def\bC{{\mathbb C}}
\def\bZ{{\mathbb Z}}

\def\bR{{\mathbb R}}
\def\bH{{\mathbb H}}

\newcommand\Symp{\operatorname{Symp}}

\newcommand\Ker{\operatorname{Ker}}

\newcommand\Diff{\operatorname{Diff}}

\catcode`\@=\active 

\begin{document}

\title{The geometry of symplectic pairs}
\author{G.~Bande}
\address{Universit\`a degli Studi di Cagliari, Dip.~Mat., Via Ospedale 
72, 09129 Cagliari, Italy}
\email{gbande{\char'100}unica.it}
\author{D.~Kotschick}
\address{Mathematisches Institut, Ludwig-Maximilians-Universit\"at M\"unchen,
Theresienstr.~39, 80333 M\"unchen, Germany}
\email{dieter{\char'100}member.ams.org}


\thanks{The authors are members of the {\sl European Differential
Geometry Endeavour} (EDGE), Research Training Network HPRN-CT-2000-00101,
supported by The European Human Potential Programme}
\date{19 July 2004; MSC 2000 classification: primary 53C15, 57R17, 57R30;
secondary 53C12, 53D35, 58A17}

\begin{abstract}
We study the geometry of manifolds carrying symplectic pairs consisting 
of two closed $2$-forms of constant ranks, whose kernel foliations are 
complementary. Using a variation of the construction of Boothby and Wang we 
build contact-symplectic and contact pairs from symplectic pairs. 
\end{abstract}

\maketitle

\section{Introduction}

A symplectic pair on a smooth manifold $M$ is a pair of non-trivial closed
two-forms $\omega_{1}$, $\omega_{2}$ of constant and complementary ranks, 
for which $\omega_{1}$ restricts as a symplectic form to the leaves of the 
kernel foliation of $\omega_{2}$, and vice versa. This definition is 
analogous to that of contact pairs and of contact-symplectic pairs 
introduced by the first author~\cite{Bande1,Bande2}. In this paper we 
exhibit several constructions of symplectic pairs on closed manifolds, and use 
them to show that even in dimension four there is a surprisingly rich supply of 
examples, with very different geometric features. One reason why this is 
suprising is that it seems that not many explicit examples are known of closed 
four-manifolds which admit pairs of complementary two-dimensional foliations. 
It was only recently that Morita and the second author constructed certain 
interesting examples within the class of foliated bundles, cf.~\cite{KM}. 
Another reason why the plethora of symplectic pairs $(\omega_{1},\omega_{2})$ 
on closed four-manifolds is surprising is that they give rise to symplectic forms 
$\omega_{\pm}=\omega_{1}\pm\omega_{2}$ compatible with different orientations 
of $M$. It is known that manifolds which are symplectic for both choices of 
orientation, or just have non-trivial Seiberg--Witten invariants for both 
orientations, are rather special, see for example~\cite{Ko}. 

On a four-manifold a symplectic pair $(\omega_{1},\omega_{2})$ can be 
equivalently defined by two symplectic forms $\omega_{\pm}$ with the 
properties $\omega_{+}\wedge\omega_{-}=0$ and 
$\omega_{+}\wedge\omega_{+}=-\omega_{-}\wedge\omega_{-}$.
If we were to change the sign in the last condition to 
$\omega_{+}\wedge\omega_{+}=\omega_{-}\wedge\omega_{-}$,
we would obtain the definition of a conformal symplectic 
couple in the sense of Geiges~\cite{Geiges2}, who gave a complete 
classification of the diffeomorphism types of four-manifolds which 
admit such couples. We will see that symplectic pairs, whose 
definition is very similar, are much more common than conformal 
symplectic couples.

While we believe that this paper shows symplectic pairs to be interesting 
geometric objects in their own right, the original motivation for this work
came from two other sources. Firstly, symplectic pairs appear naturally in 
the study of Riemannian metrics for which all products of harmonic forms 
are harmonic, see~\cite{K}, and in the investigation of the group cohomology 
of symplectomorphism groups, see~\cite{KM}. Secondly, symplectic pairs can 
be used to construct new examples of contact-symplectic and of contact pairs 
in the sense of~\cite{Bande1,Bande2,BH}. In Section~\ref{s:BW} we show how a 
variation of the classical Boothby--Wang construction~\cite{BW} allows one to 
construct contact-symplectic pairs from symplectic pairs for which $\omega_{1}$ 
represents an integral cohomology class in $M$, and contact pairs from 
contact-symplectic pairs in which the leafwise symplectic form represents an 
integral cohomology class. In particular, if we have a symplectic pair for which 
both $\omega_{i}$ represent integral classes, then the fiber product of the 
corresponding Boothby--Wang fibrations yields a contact pair on a principal 
$T^{2}$-bundle over $M$.

In Section~\ref{s:ex} we give several constructions of symplectic pairs. 
Considering their Boothby--Wang fibrations one obtains many new examples 
of contact-symplectic and of contact pairs which go beyond the examples 
exhibited in~\cite{Bande1,Bande2,BH}.
In Section~\ref{s:metric} we study Riemannian metrics compatible with 
a symplectic pair, and in Section~\ref{s:applic} we give some further 
applications of our constructions. 

In a sequel to this paper written jointly with Ghiggini~\cite{BGK}, 
we formulate and prove the appropriate adaptation of Moser's 
theorem~\cite{Moser} for symplectic pairs.

\section{The Boothby--Wang construction}\label{s:BW}

Let $(M,\omega)$ be a closed symplectic manifold. After a small 
perturbation and multiplication by a constant we may assume that 
$\omega$ represents an integral class in $H^{2}(M;\bR)$. Let 
$\pi\colon E\rightarrow M$ be the principal $S^{1}$-bundle with Euler 
class $[\omega]$. There is a connection $1$-form $\alpha$ on this 
circle bundle with curvature $\omega$, i.~e.~we have 
$d\alpha=\pi^{*}\omega$. As $\omega$ is assumed to be symplectic on $M$, 
it follows that $\alpha$ is a contact form on the total space $E$. 
This is the construction of Boothby and Wang~\cite{BW} associating so-called 
regular contact forms to integral symplectic forms.

Now if $\omega$ is an arbitrary closed $2$-form representing an 
integral cohomology class, we can again find a connection $1$-form 
$\alpha$ with curvature $\omega$, because every closed $2$-form 
representing an integral class is the curvature of some connection. 
If $\omega$ has constant rank $2k$, then it follows that $\alpha$ has 
constant class $2k+1$, meaning $\alpha\wedge (d\alpha)^{k}\neq 0$, 
and $(d\alpha)^{k+1}\equiv 0$. Recall the definition of a 
contact-symplectic pair:
\begin{definition}[\cite{Bande1,Bande2}]
    A contact-symplectic pair on a manifold $N$ consists of a $1$-form 
    $\alpha$ of constant class $2k+1$ and a closed $2$-form $\beta$ of 
    constant rank $2l$ such that the kernel foliations of $\alpha\wedge 
    (d\alpha)^{k}$ and of $\beta$ are complementary, $\alpha$ restricts as a 
    contact form to the leaves of the kernel foliation of $\beta$, and 
    $\beta$ restricts as a symplectic form to the leaves of the kernel 
    foliation of $\alpha\wedge (d\alpha)^{k}$.
    \end{definition}
Note that the kernel distributions are integrable because the forms 
$\beta$ and $\alpha\wedge (d\alpha)^{k}$ are closed. The assumption 
that the kernel foliations are complementary implies that the 
dimension of $N$ must be $2k+2l+1$.

Our discussion above immediately yields the following:
\begin{theorem}\label{t:BW1}
    Let $M$ be a closed manifold with a symplectic pair 
    $(\omega_{1},\omega_{2})$. If $[\omega_{1}]\in H^{2}(M;\bR)$ is an 
    integral cohomology class, then the total space of the circle 
    bundle $\pi\colon E\rightarrow M$ with Euler class $[\omega_{1}]$ 
    carries a natural contact-symplectic pair.
    \end{theorem}
Indeed, if $\alpha$ is a connection form with curvature $\omega_{1}$
and $\beta=\pi^{*}\omega_{2}$, then all the required properties are 
satisfied.

Note that if ${\mathcal F_{i}}$ denotes the kernel foliation of 
$\omega_{i}$, then in the above construction the classical 
Boothby--Wang construction is performed leafwise over the leaves of 
${\mathcal F_{2}}$, to which $\omega_{1}$ restricts as an integral 
symplectic form. On $E$ the kernel foliation of 
$\beta=\pi^{*}\omega_{2}$ consists of the circle bundles over the 
leaves of ${\mathcal F_{2}}$, whereas the kernel foliation of 
$\alpha\wedge (d\alpha)^{k}$ is complementary and obtained by lifting 
the leaves of ${\mathcal F_{1}}$ to the horizontal subspaces for the 
connection $\alpha$.

Next recall the definition of a contact pair:
\begin{definition}[\cite{Bande1,BH}]
    A contact pair on a manifold $N$ consists of a pair of $1$-forms 
    $(\alpha,\gamma)$ of constant class $2k+1$ and $2l+1$ 
    respectively, such that the kernel foliations of $\alpha\wedge 
    (d\alpha)^{k}$ and of $\gamma\wedge (d\gamma)^{l}$ are complementary, 
    $\alpha$ restricts as a contact form to the leaves of the kernel 
    foliation of $\gamma\wedge (d\gamma)^{l}$, and $\gamma$ restricts as a 
    contact form to the leaves of the kernel foliation of $\alpha\wedge 
    (d\alpha)^{k}$.
    \end{definition}
Again the kernel distributions are integrable because the defining forms 
are closed. The assumption that the kernel foliations are complementary 
implies that the dimension of $N$ must be $2k+2l+2$.

Our discussion above yields the following:
\begin{theorem}\label{t:BW2}
    Let $M$ be a closed manifold with a contact-symplectic pair 
    $(\alpha,\beta)$. If $[\beta]\in H^{2}(M;\bR)$ is an integral 
    cohomology class, then the total space of the circle bundle 
    $\pi\colon E\rightarrow M$ with Euler class $[\beta]$ 
    carries a natural contact pair.
    \end{theorem}
Indeed, if $\gamma$ is a connection form with curvature $\beta$
and we identify $\alpha$ with its pullback under $\pi$, then all 
the required properties are satisfied.

Finally, combining Theorems~\ref{t:BW1} and~\ref{t:BW2}, we obtain:
\begin{corollary}
    If a closed manifold $M$ has a symplectic pair 
    $(\omega_{1},\omega_{2})$ such that both $[\omega_{i}]\in 
    H^{2}(M;\bR)$ are integral, then the fiber product of the two 
    circle bundles with Euler classes equal to $[\omega_{1}]$ and 
    $[\omega_{2}]$ respectively carries a natural contact pair.
    \end{corollary}

An important difference between this leafwise Boothby--Wang 
construction and the classical one is that we cannot perturb the 
defining forms in a symplectic pair so as to make them rational\footnote{and 
integral after multiplication with a constant}, because one cannot 
control the rank under such perturbations. Therefore, we will check in 
each example of a symplectic pair we construct in the next section, 
whether the defining forms represent integral cohomology classes. 

\section{Constructions of symplectic pairs}\label{s:ex}

The most obvious examples of symplectic pairs are of course products 
of symplectic manifolds with the induced split symplectic forms. In 
this case one can obviously choose the forms $\omega_{i}$ to represent 
integral classes.

We now discuss non-trival sources of examples.

\subsection{Flat bundles with symplectic total holonomy}
Let $(B,\omega_{B})$ and $(F,\omega_{F})$ be closed symplectic 
manifolds, and $\rho\colon\pi_{1}(B)\rightarrow\Symp(F,\omega_{F})$ a 
representation of the fundamental group of $B$ in the group of
symplectomorphisms of $(F,\omega_{F})$. The suspension of $\rho$ 
defines a horizontal foliation on the fiber bundle $\pi\colon 
M\rightarrow B$ with fiber $F$ and total space $M=(\tilde 
B\times F)/\pi_{1}(B)$, where $\pi_{1}(B)$ acts on $\tilde B$ by 
covering transformations and on $F$ via $\rho$. As the image of $\rho$ 
preserves the symplectic form $\omega_{F}$, the pullback of this form 
to the product $\tilde B\times F$ descends to $M$ as a closed form of 
constant rank, whose kernel foliation is exactly the horizontal 
foliation complementary to the fibers. Pulling back $\omega_{B}$ to 
the total space $M$ we obtain another closed form of constant rank, 
which is a defining form for the vertical foliation whose leaves are the 
fibers of the fibration. As the two foliations are complementary by 
construction, the forms $\omega_{F}$ and $\pi^{*}\omega_{B}$ form a 
symplectic pair on $M$.

Note that if we choose $\omega_{B}$ to be integral, then so is its 
pullback. For $\omega_{F}$ checking integrality is more subtle. In 
particular, it turns out that starting with an integral form on $F$, 
though necessary, is not usually sufficient. 

A special case of the above construction is given by taking a single 
symplectic diffeomorphism $\varphi\in\Symp(F,\omega_{F})$, and forming 
the product of its mapping torus $M_{\varphi}$ with $S^{1}$. 
If $\varphi$ is isotopic to the identity through symplectomorphisms 
$\varphi_{t}$, with $\varphi_{1}=\varphi$ and $\varphi_{0}=Id_{F}$, 
then $M_{\varphi}$ is diffeomorphic to $F\times S^{1}$ by a 
diffeomorphism encoding the isotopy. It was proved in Lemma 8 
of~\cite{KM} that under this diffeomorphism, the cohomology class 
$[\omega_{F}]\in H^{2}(M_{\varphi})$ corresponds to 
$[\omega_{F}]+Flux(\varphi_{t})\otimes\nu\in H^{2}(F)\oplus 
(H^{1}(F)\otimes H^{1}(S^{1}))$, where $\nu$ is the fundamental 
cohomology class of $S^{1}$. Thus, the cohomology class of 
$\omega_{F}$ on such a symplectic mapping torus is integral if 
and only if $Flux(\varphi_{t})$ is an integral class in $H^{1}(F)$.

While symplectic mapping tori have a rather simple topology, 
determined completely by $\varphi$, there are more complicated flat 
bundles with symplectic total holonomy which exhibit more complex 
topology. For example, in the simplest possible case, where $B$ and 
$F$ are both $2$-dimensional, Kotschick and Morita~\cite{KM} proved 
the following:
\begin{theorem}[\cite{KM}]\label{t:KM}
    For every $g\geq 3$ there exist foliated oriented surface bundles
$\pi\colon M \ra B$ over closed oriented surfaces $B$ with fibers $F$ 
of genus $g$, which have non-zero signature and whose total holonomy 
group is contained in the symplectomorphism group $\Symp(F,\omega_{F})$ 
with respect to a prescribed symplectic form $\omega_{F}$ on $F$. In 
fact, one can restrict the holonomy to be in $\Symp(F;D^{2})$, the 
group of compactly supported symplectomorphisms of $F\setminus D^{2}$.
    \end{theorem}
The first part is Theorem~1 in~\cite{KM}, whereas the addendum restricting 
to symplectomorphisms relative to an embedded disk follows from the proof 
of Theorem~3 in~\cite{KM}. This addendum is useful for the construction of 
further symplectic pairs, see~\ref{ss:sums} below, because it 
implies that the $4$-manifold $M$ in the statement of the theorem 
contains a product neighbourhood $D^{2}\times B$ to which the 
symplectic pairs restrict in the obvious way, so that the two 
foliations are given by $D^{2}\times\{\star\}$ and by $\{\star\}\times B$. 
In particular, the horizontal foliation has an open set of closed 
leaves.

\begin{remark}
    If the base manifold $B$ is not just symplectic, but has a 
    symplectic pair, then any flat bundle over $B$ with symplectic 
    total holonomy inherits something we may naturally call a symplectic 
    triple. From this one can combine several different symplectic pairs. 
    The same remark applies if the total holonomy preserves a symplectic 
    pair on the fiber $F$. 
    \end{remark}

Foliated bundles can also be used to construct contact-symplectic and 
contact pairs directly. For example, if $B$ carries a contact or symplectic 
structure and the image of a homomorphism $\rho\colon\pi_{1}(B)\rightarrow\Diff(F)$ 
preserves a contact form on $F$, then $M=(\tilde B\times F)/\pi_{1}(B)$ 
obtained by suspending $\rho$ inherits a contact or contact-symplectic pair.

\subsection{The Gompf sum for symplectic pairs}\label{ss:sums}

Gompf~\cite{Gompf} has shown that symplectic manifolds with closed 
symplectic submanifolds of codimension $2$ admit certain cut-and-paste 
constructions which build new symplectic manifolds out of old ones. 
Suppose that $(M_{1},\omega_{1})$ and $(M_{2},\omega_{2})$ are closed 
symplectic manifolds of dimension $2n$ admitting symplectic 
submanifolds $\Sigma_{i}\subset M_{i}$ of codimension $2$ with 
trivial normal bundles, and such that $(\Sigma_{1},\omega_{1})$ and 
$(\Sigma_{2},\omega_{2})$ are symplectomorphic. Then by the symplectic 
tubular neighbourhood theorem they have symplectomorphic neighbourhoods. 
In this situation $M_{1}\setminus\Sigma_{1}$ and $M_{2}\setminus\Sigma_{2}$ can 
be glued together symplectically along punctured tubular neighbourhoods 
of the $\Sigma_{i}$. The gluing map turns a punctured normal disk inside 
out symplectically. 

This construction sometimes works for manifolds with symplectic pairs 
if one of the foliations has codimension $2$ and has an open set of compact 
leaves. Let $(M_{1},\omega_{1},\omega_{2})$ and $(M_{2},\eta_{1},\eta_{2})$ 
be closed manifolds of dimension $2n$ with symplectic pairs for which 
$rk(\omega_{1})=rk(\eta_{1})=2n-2$. Suppose that the kernel foliations 
$\FF_{1}$ of $\omega_{2}$ and $\FF_{2}$ of $\eta_{2}$ each have closed leaves 
$\Sigma_{1}$ and $\Sigma_{2}$ respectively, such that our symplectic pairs 
admit product structures in open neighbourhoods of the $\Sigma_{i}$. 
This means that we assume that $\Sigma_{1}$ has an open neighbourhood 
$U_{1}\subset M_{1}$ which is diffeomorphic to $\Sigma_{1}\times D^{2}$ 
in such a way that 
$\omega_{1}\vert_{U_{1}}=\pi_{1}^{*}(\omega_{1}\vert_{\Sigma_{1}})$ 
and $\omega_{2}\vert_{U_{1}}=\pi_{2}^{*}(\omega_{2}\vert_{D^{2}})$, 
where the $\pi_{i}$ are the projections to the factors; and similarly 
for $\Sigma_{2}\subset M_{2}$. Then we may assume without further loss 
of generality that $\omega_{2}\vert_{D^{2}}$ and $\eta_{2}\vert_{D^{2}}$ 
coincide with the standard area form $dx\wedge dy$ on the disk.
Suppose further that there is a symplectomorphism 
    $$
    f\colon (\Sigma_{1},\omega_{1})\longrightarrow 
    (\Sigma_{2},\eta_{1}) \ .
    $$
Then the Gompf sum $M_{1}\natural_{f}M_{2}$ of the $M_{i}$ along 
the submanifolds $\Sigma_{i}$ carries a natural symplectic pair.

As in Gompf's original construction~\cite{Gompf}, the assumptions are 
particularly easy to verify when the $M_{i}$ are $4$-dimensional. In 
this case $\omega_{1}$ and $\eta_{1}$ are volume forms on the $\Sigma_{i}$ 
and, by Moser's theorem~\cite{Moser}, a symplectomorphism $f$ as 
above exists as soon as $\Sigma_{1}$ and $\Sigma_{2}$ have the same 
genus and 
$$
\int_{\Sigma_{1}}\omega_{1} = \int_{\Sigma_{2}}\eta_{1} \ .
$$

We can use the flat bundles in Theorem~\ref{t:KM} as building blocks for 
the Gompf sum, because, by construction, their horizontal foliations have 
product structures on an open set. For the vertical foliations we trivially 
have product structures around every fiber. Performing the Gompf sum of 
symplectic pairs by matching fibers with fibers or sections\footnote{ $=$ 
closed leaves of the horizontal foliations} with sections does not lead to 
any new examples. However, taking a flat bundle over a surface of genus $g$, 
and another one with fibers of genus $g$, we can, after scaling one of the 
$2$-forms involved by a constant, perform the sum of symplectic pairs matching 
a fiber in one fibration with the section in the other fibration. This gives 
new examples of manifolds admitting symplectic pairs which are not surface 
bundles over surfaces. 

\subsection{Four-dimensional Thurston geometries}

A geometry in the sense of Thurston consists of a model space $X$ 
which is a simply connected complete Riemannian manifold, together with 
a group $G$ of effective isometries acting transitively and admitting a 
discrete subgroup $\Gamma$ for which the quotient space $X/\Gamma$ is a 
compact smooth manifold. Such compact quotients are said to admit a 
Thurston geometry of type $(X,G)$.

The four-dimensional Thurston geometries have been classified by 
Filipkiewicz (unpublished). We refer the reader to Wall's 
papers~\cite{W1,W2} for an account of this classification. We now want to 
show that for some of these geometries there are natural $G$-invariant 
symplectic pairs on the model spaces, which then descend to all compact 
quotients. As the isometries we consider preserve a symplectic structure, 
they are orientation-preserving.

\begin{example}\label{ex:prod1}
    Consider the model spaces $S^{2}\times\bR^{2}$, $S^{2}\times\bH^{2}$, 
    $\bR^{2}\times\bH^{2}$ with the product metrics obtained from the 
    standard constant curvature metrics on the factors. In this case any 
    isometry preserves the local product structure, and its factors. 
    In the maximal group of orientation-preserving isometries those 
    which preserve a pair of given orientations on the factors form a 
    subgroup $G$ of index $2$. The volume forms of the metrics on the 
    factors then form a $G$-invariant symplectic pair on $X$.
    \end{example}

\begin{example}\label{ex:prod2}
    The discussion in the previous example applies to the model spaces 
    $S^{2}\times S^{2}$ and $\bH^{2}\times\bH^{2}$, except that these 
    also admit isometries interchanging the two factors.
    \end{example}

\begin{example}\label{ex:flat}
    The model space $\bR^{4}$ with its standard flat metric has as its 
    compact quotients the flat Riemannian $4$-manifolds $M$. If such 
    a manifold is orientable, then $b_{1}(M)>0$. It is known that 
    $b_{1}(M)\leq 4$, with equality if and only if $M$ is 
    diffeomorphic to $T^{4}$, and that $b_{1}(M)\neq 3$. Moreover, 
    if $b_{1}(M)=1$, then the vanishing of the Euler characteristic 
    shows that $b_{2}(M)=0$, so that $M$ cannot be symplectic. Thus 
    the only interesting case is when $b_{1}(M)=2$. The classification 
    of flat $4$-manifolds in~\cite{Hil,wagner} shows that in the case 
    $b_{1}(M)=2$ they are all quotients of $\bR^{4}$ by isometry groups 
    preserving a product structure $\bR^{2}\times\bR^{2}$, and acting 
    on each factor preserving its orientation. In a different guise, 
    this statement appears in the classification of compact complex 
    surfaces, where these particular flat Riemannian manifolds appear 
    as so-called hyperelliptic\footnote{These surfaces are sometimes 
    called bielliptic, because they have two different elliptic 
    fibrations.} sufaces, see~\cite{BPV} p.~148. They are in fact quotients 
    of products of elliptic curves by free diagonal actions of finite groups 
    of holomorphic automorphisms. Thus they carry natural symplectic pairs. 
    \end{example}

\begin{example}\label{ex:sol}
    Consider the model space $X=Sol^{3}\times\bR$ with its maximally 
    symmetric product metric. Then the maximal connected isometry 
    group $G_{0}$ is also $Sol^{3}\times\bR$, acting on itself by left 
    multiplication, cf.~\cite{Ue} p.~518/19. This Lie group admits a 
    parallelization by left-invariant one-forms 
    $\alpha_{1},\ldots,\alpha_{4}$ with 
    $d\alpha_{1}=\alpha_{1}\wedge\alpha_{4}$, 
    $d\alpha_{3}=\alpha_{4}\wedge\alpha_{3}$, 
    $d\alpha_{2}=d\alpha_{4}=0$. It follows that 
    $\omega_{1}=\alpha_{1}\wedge\alpha_{3}$ and 
    $\omega_{2}=\alpha_{2}\wedge\alpha_{4}$ form a left-invariant 
    symplectic pair.
    \end{example}
 
\begin{example}\label{ex:nil4}
    Consider the model space $X=Nil^{4}$ with its maximally 
    symmetric metric. Then again the maximal connected isometry 
    group $G_{0}$ coincides with $X$, acting on itself by left 
    multiplication, cf.~\cite{Ue} p.~518. This Lie group admits a 
    parallelization by left-invariant one-forms 
    $\alpha_{1},\ldots,\alpha_{4}$ with 
    $d\alpha_{2}=\alpha_{1}\wedge\alpha_{4}$, 
    $d\alpha_{3}=\alpha_{2}\wedge\alpha_{4}$, 
    $d\alpha_{1}=d\alpha_{4}=0$. It follows that 
    $\omega_{1}=\alpha_{1}\wedge\alpha_{2}$ and 
    $\omega_{2}=\alpha_{3}\wedge\alpha_{4}$ form a left-invariant 
    symplectic pair.
    \end{example}

\begin{example}\label{ex:nil3}
    Consider the model space $X=Nil^{3}\times\bR$ with its maximally 
    symmetric product metric. This Lie group admits a parallelization by 
    left-invariant one-forms $\alpha_{1},\ldots,\alpha_{4}$ with 
    $d\alpha_{3}=\alpha_{1}\wedge\alpha_{2}$ and
    $d\alpha_{1}=d\alpha_{2}=d\alpha_{4}=0$. It follows that 
    $\omega_{1}=\alpha_{1}\wedge\alpha_{3}$ and 
    $\omega_{2}=\alpha_{2}\wedge\alpha_{4}$ form a left-invariant 
    symplectic pair. In this case the maximal connected group of isometries 
    is larger than $Nil^{3}\times\bR$, because it contains the 
    rotations in the plane spanned by $\alpha_{1}$ and $\alpha_{2}$. 
    But these rotations do not preserve the symplectic pair.
    \end{example}

It turns out that the remaining Thurston geometries do not support 
any symplectic pairs:
\begin{theorem}
    The model spaces $S^{4}$, $\bC P^{2}$, $\bH^{4}$, $\bC H^{2}$, 
    $\widetilde{PSL_{2}(\bR)}\times\bR$, $\bH^{3}\times\bR$, 
    $S^{3}\times\bR$, $Sol^{4}_{0}$, $Sol^{4}_{1}$ and $Sol_{m,n}$ 
    with $m\neq n$ with their standard metrics do not admit any 
    transitive groups of isometries containing cocompact lattices 
    which also preserve a symplectic pair.
    \end{theorem}
We have formulated the theorem in such a way that it covers 
non-maximal geometries in the sense of~\cite{W1,W2}, i.~e.~we rule 
out symplectic pairs invariant under transitive subgroups which need 
not be the maximal isometry groups.
\begin{proof}
    We proceed case-by-case. The four-sphere admits no symplectic 
    structure, and so is ruled out. Any compact quotient of 
    $S^{3}\times\bR$ is finitely covered by $S^{3}\times S^{1}$, and 
    so admits no symplectic structure.
    
    The complex projective plane does 
    admit a symplectic structure, but its tangent bundle has no 
    decomposition into a direct sum of two oriented plane bundles. 
    (This is equivalent to the well-known fact that $\bC P^{2}$ 
    endowed with the non-complex orientation admits no almost-complex 
    structure.) Thus $\bC P^{2}$ is also ruled out. Concerning its 
    non-compact dual $\bC H^{2}$, Wall~\cite{W2} proved that the 
    isotropy subgroup of any transitive isometry group admitting a 
    cocompact lattice contains $U(2)$. As this does not preserve any 
    splitting of $\R^{4}$ into a direct sum of proper subspaces, $\bC 
    H^{2}$ cannot carry any invariant symplectic pair.
    
    This last argument also applies to the geometries $\bH^{4}$ and 
    $\bH^{3}\times\bR$. In these cases the isotropy subgroup of any 
    transitive isometry group admitting a cocompact lattice contains 
    $SO(4)$ respectively $SO(3)$. These groups do not preserve any 
    splitting of $\R^{4}$ into a direct sum of $2$-dimensional subspaces.

    Finally, for the Lie group geometries $\widetilde{PSL_{2}(\bR)}\times\bR$, 
    $Sol^{4}_{0}$, $Sol^{4}_{1}$ and $Sol_{m,n}$ with $m\neq n$, 
    any transitive isometry group must contain the Lie group itself, 
    acting by left multiplication. However, in these cases it is easy 
    to check using the structure constants in~\cite{W1} that there are 
    no left-invariant symplectic forms, cf.~\cite{Hamilton}.
    \end{proof}
 
\section{Compatible metrics}\label{s:metric}

In this section we clarify the metric properties of symplectic pairs. 
As a first step, we have the following:
\begin{proposition}
    Let $M$ be a manifold endowed with two smooth complementary 
    foliations $\FF$ and $\GG$ which admit closed defining forms. 
    Then there are Riemannian metrics $g$ on $M$ for which $\FF$ 
    and $\GG$ are orthogonal and have minimal leaves.
    \end{proposition}
\begin{proof}
    This is a consequence of the minimality criterion of Rummler and 
    Sullivan, see~\cite{God}, page 371/372. Given an arbitrary 
    foliation $\FF$ with leaves of dimension $d$ and a form of degree 
    $d$ which is relatively closed for $\FF$ and restricts as a 
    volume form to the leaves of $\FF$, one can construct metrics $g$ 
    making the leaves of $\FF$ minimal, and such that the given 
    $d$-form is the volume form of the restricted metric. These 
    metrics $g$ can be chosen to make the kernel of the $d$-form 
    orthogonal to $\FF$, and the restriction to this orthogonal 
    complement is arbitrary.
    
    Suppose that $\FF=\Ker(\alpha)$ and $\GG=\Ker(\beta)$, with $\alpha$ 
    and $\beta$ closed and of degrees equal to the codimensions of $\FF$ 
    and $\GG$ respectively. As $\FF$ and $\GG$ are assumed to be 
    complementary, $\alpha$ is a leafwise volume form on $\GG$ and 
    $\beta$ is a leafwise volume form on $\FF$. Define a metric $g$ by 
    requiring $T\FF$ and $T\GG$ to be orthogonal, and choosing $g$ 
    along $\FF$ so that $\beta$ is the Riemannian volume form of 
    $g\vert T\FF$, and choosing $g$ along $\GG$ so that $\alpha$ is 
    the Riemannian volume form of $g\vert T\GG$. These requirements 
    clearly underdetermine the metric, and any such metric has all the 
    desired properties.
    \end{proof}

\begin{corollary}\label{cor}
    A manifold endowed with a symplectic, contact-symplectic or 
    contact pair admits metrics for which the characteristic 
    foliations are orthogonal with minimal leaves.
    \end{corollary}
    
In many of the examples constructed above there are metrics which in 
addition to making the foliations orthogonal with minimal leaves have 
further good properties. For example, the flat bundles 
always have metrics for which the vertical foliation 
is Riemannian and the horizontal foliation has totally geodesic leaves. 
The following theorem shows that a general symplectic pair does not 
admit any metric with properties more restrictive than the ones 
specified in Corollary~\ref{cor}.
\begin{theorem}
    There are symplectic pairs on closed four-manifolds for which 
    both foliations are not geodesible and not Riemannian.
    \end{theorem}
\begin{proof}
    Consider foliated surface bundles $M$ over surfaces with 
    symplectic total holonomy. The normal bundle of the horizontal 
    foliation is the tangent bundle along the fibers, and its first
    Pontryagin number is three times the signature $\sigma(M)$, because 
    the Pontryagin number of the tangent bundle of the horizontal 
    foliation vanishes. If the signature is non-zero, then Pasternack's 
    refinement~\cite{P} of the Bott vanishing theorem for Riemannian foliations 
    implies that the horizontal foliation is not Riemannian. To see 
    this, recall that for Riemannian foliations Pasternack shows that 
    the Pontryagin numbers of the normal bundle vanish in degrees 
    strictly larger than the codimension of the foliation, which 
    improves the range of vanishing in Bott's theorem by a factor of 
    two. In our situation this means that the first Pontryagin number 
    of the normal bundle, which is in degree $4$, vanishes, as the 
    codimension equals $2$.
    
    Now take two such foliated bundles, $M_{1}$ and $M_{2}$. By 
    Theorem~\ref{t:KM} we can choose both of them with non-zero 
    signature, such that the base genus of $M_{2}$ equals the 
    fiber genus of $M_{1}$, and such that the horizontal foliation 
    in $M_{2}$ has an open set of compact leaves with trivial normal 
    bundle. Let $M$ be the Gompf sum $M_{1}\natural_{f}M_{2}$, where 
    a section in $M_{2}$ is identified with a fiber in $M_{1}$, as 
    discussed in Subsection~\ref{ss:sums} above. This sum $M$ carries 
    an induced symplectic pair, and it is clear that on $M$ the first 
    Pontryagin number of both $T\FF$ and $T\GG$ is non-zero, because 
    we have chosen both $M_{i}$ to have non-zero signature. Note that 
    each of these bundles is the normal bundle for the complementary 
    foliation. Thus neither of the two foliations can be Riemannian, 
    by Pasternack's theorem~\cite{P}.
    
    Suppose now that in $M$ one of the foliations, say $\FF$, is 
    geodesible. If a metric making $\FF$ totally geodesic also makes 
    it orthogonal to $\GG$, then the duality theorem for totally 
    geodesic and bundle-like foliations implies that $\GG$ is 
    Riemannian, see~\cite{God} p.~190. This is a contradiction. 
    
    Next assume that we can choose a metric $g$ for which $\FF$ is 
    totally geodesic, without assuming that its orthogonal complement 
    is $\GG$. Cairns and Ghys~\cite{CG} have shown that for any 
    two-dimensional geodesible foliation on a $4$-manifold we may 
    choose $g$ to make the leaves both totally geodesic and of constant 
    Gaussian curvature. As $\FF$ has closed leaves of genus $\geq 2$, the constant 
    curvature is negative. Another result of~\cite{CG} then tells us 
    that the $g$-orthogonal complement $T\FF^{\perp}$ is integrable, 
    and defines a foliation $\HHH$ (which may be different from $\GG$). 
    By the duality theorem, $g$ is bundle-like for $\HHH$. 
    But $\HHH$ has normal bundle $T\FF$, which has non-zero first 
    Pontryagin number, and so we again have a contradiction with 
    Pasternack's theorem.
\end{proof}
There are special cases of symplectic pairs for which it is possible 
to find a metric which makes the two foliations orthogonal and 
totally geodesic, for example the Thurston geometries which are 
products of two-dimensional geometries. When performing a 
Boothby--Wang construction on such an example one can choose a 
submersion metric on the total space which also has the property that 
the foliations of the contact-symplectic pair are orthogonal and 
totally geodesic. This will be used in Subsection~\ref{ss:ex} below.

\section{Some applications}\label{s:applic}

\subsection{Torus-bundles over the torus}

We now want to prove the following:
\begin{theorem}
    Every oriented $T^{2}$-bundle over $T^{2}$ admits a symplectic 
    pair $(\omega_{1},\omega_{2})$ for which the cohomology classes 
    of the $\omega_{i}$ are integral.
    \end{theorem}
This can be seen as generalizing a result of Geiges~\cite{Geiges}, 
who proved that these manifolds admit symplectic structures. His proof, 
like ours, depends in the classification of $T^{2}$-bundles over $T^{2}$ 
due to Sakamoto and Fukuhara~\cite{SF}, and on the fact that all these 
manifolds carry compatible Thurston geometries, cf.~\cite{Ue}.
\begin{proof}
    The classification of orientable $T^{2}$-bundles over $T^{2}$ is 
    summarized in the table in the appendix. We will proceed 
    case-by-case and use the information given in the table.
    In case (a), for the four-torus, the claim is trivial.
    
    Case (b) consists of manifolds with Thurston geometry 
    $Nil^{3}\times\bR$. As the first Betti number equals $3$, these 
    manifolds are nilmanifolds (rather than infranil manifolds), 
    i.~e.~they are quotients of our Lie group by lattices in the 
    group itself acting by left translations, cf.~\cite{H} p.~170.
    We saw in Example~\ref{ex:nil3} that there is a left-invariant 
    symplectic pair on the group. Thus this descends to all manifolds 
    under discussion here.
    
    In case (c) we have the flat orientable four-manifolds with 
    $b_{1}=2$. These have symplectic pairs by Example~\ref{ex:flat}.
    
    Case (d) consists of manifolds with Thurston geometry $Nil^{4}$. 
    As their first Betti number equals $2$, these manifolds are nilmanifolds 
    (rather than infranil manifolds), cf.~\cite{H} p.~170.
    We saw in Example~\ref{ex:nil4} that there is a left-invariant 
    symplectic pair on the group. Thus this descends to all manifolds 
    under discussion in this case.
    
    Cases (e) and (f) consist of infranil manifolds for the group 
    $Nil^{3}\times\bR$. As we do not have a symplectic pair on the 
    model space invariant under the full group of orientation-preserving 
    isometries, we argue instead as in~\cite{Geiges}. Geiges~\cite{Geiges} 
    showed that identifying the model space with $\R^{4}$ with coordinates 
    $(x,y,z,t)$, the two-forms $dy\wedge dt$ and $dx\wedge dz - xdx\wedge dy$ 
    are invariant under the lattices arising as fundamental groups in this 
    case. Clearly they are closed of constant rank equal to $2$, and their wedge 
    product is a volume form. Thus they give rise to a symplectic pair.
    
    Cases (g) and (h) consist of manifolds with Thurston geometry 
    $Sol^{3}\times\bR$. It was shown in~\cite{Geiges} that identifying 
    the model space with $\R^{4}$ with coordinates $(x,y,z,t)$, the two-forms 
    $dx\wedge dy$ and $dz\wedge dt$ are invariant under the lattices arising 
    as fundamental groups. They are closed of constant rank equal to $2$, and 
    their wedge product is a volume form. Thus they give rise to a symplectic pair.
    
It remains to address the integrality of the cohomology classes of the forms 
involved. This can trivially be arranged in the case of $T^{4}$. For the 
nilmanifolds of $Nil^{3}\times\bR$ the integrality of the cohomology classes 
for the symplectic pair we have exhibited can be checked by direct calculation, 
or using~\cite{Geiges}. For the remaining cases, (c) -- (h), we give a uniform 
argument as follows. All $T^{2}$-bundles over $T^{2}$ have vanishing  Euler 
characteristic and signature. Therefore, if $b_{1}=2$, we conclude that $b_{2}=2$, 
and the intersection form is indefinite. Thus, $H^{2}(M;\bR)$ equipped with the 
cup product form is hyperbolic, and the classes of square $0$ make up the light 
cone. For a symplectic pair the classes $[\omega_{i}]$ have square zero and 
$[\omega_{1}]\cdot [\omega_{2}]\neq 0$, thus they span the light cone. 
It follows that after constant rescaling these classes are integral.
\end{proof}
    
\subsection{Irreducible quotients of the polydisk}\label{ss:ex}

For our final application we return to the Thurston geometry with model 
space $\bH^{2}\times\bH^{2}$, which we discussed briefly in 
Example~\ref{ex:prod2}. The connected component of the identity in the 
isometry group is $PSL_{2}(\bR)\times PSL_{2}(\bR)$, acting on the model 
space preserving the symplectic pair formed by the volume forms $\omega_{1}$ 
and $\omega_{2}$ of the hyperbolic metrics on the factors. Note that 
the product metric on $\bH^{2}\times\bH^{2}$ is K\"ahler for both 
choices of orientation, with K\"ahler forms $\omega_{1}\pm\omega_{2}$.

It is well known that there are irreducible cocompact lattices 
$\Gamma\subset PSL_{2}(\bR)\times PSL_{2}(\bR)$, where by 
irreducibility we mean that $\Gamma$ is not commensurate to a product 
of lattices in $PSL_{2}(\bR)$. While the existence of irreducible 
lattices can be deduced from a general theorem due to Borel, there 
are actually explicit constructions in this case due to Kuga (cf.~\cite{BPV}) 
and Shavel~\cite{Sh} using the theory of quadratic forms. The quotients 
$(\bH^{2}\times\bH^{2})/\Gamma$ are compact complex-algebraic surfaces of general 
type with K\"ahler class $\omega_{1}\pm\omega_{2}$ (up to scale, the sign 
depending on the choice of orientations). It follows that the 
$\omega_{i}$ represent integral classes in cohomology.

These irreducible quotients of the polydisk have already been used 
to exhibit various interesting phenomena in both differential and 
algebraic geometry, cf.~\cite{CG,SB}. Here we shall add one more, in 
the form of the following result:
\begin{proposition}\label{p:M}
    There exist closed $5$-manifolds $M$ with two complementary 
    foliations and a Riemannian metric for which the foliations are 
    orthogonal and totally geodesic, and such that $M$ does not admit 
    any finite cover by a product of manifolds of strictly smaller 
    dimension.
    \end{proposition}
\begin{proof}
    Let $X$ be the quotient of the polydisk by a torsion-free irreducible 
    cocompact lattice $\Gamma$. This carries two complementary 
    foliations which are orthogonal and totally geodesic with respect 
    to the metric induced from the product metric on the universal 
    covering.
    
    As the cohomology class of the form $\omega_{1}$ is integral, we can 
    perform the leafwise Boothby--Wang construction of Section~\ref{s:BW} 
    to obtain a closed manifold $M$, which is the total space of the 
    corresponding circle bundle over $X$. On $M$ we obtain a contact-symplectic 
    pair, and a metric for which the two foliations are orthogonal and totally 
    geodesic. In fact, the Riemannian universal covering of $M$ is isometric 
    to the direct product $\bH^{2}\times\widetilde{PSL_{2}(\bR)}$, where we 
    think of $\widetilde{PSL_{2}(\bR)}$ as the universal covering of the unit 
    tangent bundle of $\bH^{2}$.
    
    It remains to prove that $M$ does not have any finite covering which 
    splits as a direct product of two manifolds of positive dimension. 
    Now it is known that $X$ has vanishing first Betti number, see~\cite{Sh}, 
    and therefore $M$ also has vanishing first Betti number by the 
    Gysin sequence of the circle fibration. It is easy to 
    see that the same conclusion must hold for any finite covering of $M$. 
    Thus, no such covering can split off a circle, and if it is homotopy 
    equivalent to a product of a $2$-manifold and a $3$-manifold, then these 
    factors must be real homology spheres. By the classification of 
    surfaces the $2$-dimensional factor is then $S^{2}$, contradicting the 
    fact that $M$ and all its finite coverings are aspherical.
    \end{proof}

\begin{remark}
    Note that we have excluded all splittings of finite coverings of 
    $M$, without assuming that they are induced by the foliations.
    \end{remark}
    
\begin{remark}
    Proposition~\ref{p:M} answers a question of Matveev~\cite{Matveev}, 
    related to his work in~\cite{M}. He noted that in dimensions $2$ 
    and $3$ every closed Riemannian manifold with a local product structure 
    given by a pair of orthogonal totally geodesic foliations admits a 
    finite covering which is a genuine product (not necessarily 
    induced by the foliations). In dimension $4$ this result is false 
    because of the existence of irreducible quotients of the polydisk, 
    and similar examples exist also in dimensions $\geq 6$.
    \end{remark}

\medskip
\noindent
{\bf Acknowledgement:} We are grateful to N.~A'Campo for having pointed 
out the overlap in our interests, which led to the present collaboration.

\newpage

\section*{Appendix: Orientable $T^{2}$-bundles over $T^{2}$}

The following table summarizes the classification of orientable 
$T^{2}$-bundles over $T^{2}$ due to Sakamoto and Fukuhara~\cite{SF}, 
and the information about their Thurston geometries due to 
Ue~\cite{Ue}, compare also~\cite{Geiges}.

\bigskip

{
\begin{tabular}{|c|c|l|c|} \hline
& & &  \\ 
& \raisebox{1.5ex}[-1.5ex]{$b_1$} & 
\raisebox{1.5ex}[-1.5ex]{Monodromy $\&$ Euler class} & 
\raisebox{1.5ex}[-1.5ex]{Geometry} \\ \hline\hline
{\rule [-3mm]{0mm}{8mm}(a) }& 4 & $\{I,I,(0,0)\} = T^4$ & $\mathbb{R}^4$  \\ \hline
{\rule [-3mm]{0mm}{8mm}(b) }& 3 & $\{I,I,(m,n)\}, (m,n) \neq (0,0)$ & 
$Nil^3 \times \mathbb{R}$  \\ \hline
{\rule [-3mm]{0mm}{8mm}(c) }& 2 & $\left\{ \left(\begin{array}{lr}
                                                    0 & -1 \\
                                                    1 & -1 \\
                                           \end{array} \right), I, (0,0) \right\}$ & $\mathbb{R}^4$ \\ 
{\rule [-3mm]{0mm}{8mm}}& & $\left\{ \left(\begin{array}{lr}
                                                    0 & -1 \\
                                                    1 & -1 \\
                                           \end{array} \right), I, (-1,0) \right\} $&  \\
{\rule [-3mm]{0mm}{8mm}}& & $\left\{ \left(\begin{array}{lr}
                                                    0 & -1 \\
                                                    1 & 0 \\
                                           \end{array} \right), I, (0,0) \right\} $&  \\                                           
{\rule [-3mm]{0mm}{8mm}}& & $\left\{ \left(\begin{array}{lr}
                                                    0 & -1 \\
                                                    1 & 0 \\
                                           \end{array} \right), I, (-1,0) \right\}$ &  \\
{\rule [-3mm]{0mm}{8mm}}& & $\left\{ \left(\begin{array}{lr}
                                                    1 & -1 \\
                                                    1 & 0 \\
                                           \end{array} \right), I, (0,0) \right\}$ &  \\
{\rule [-3mm]{0mm}{8mm}}& & $\{-I, I, (0,0) \}$ &  \\
{\rule [-3mm]{0mm}{8mm}}& & $\{-I, I ,(-1,0) \}$ &  \\ \hline
                                          
{\rule [-3mm]{0mm}{8mm}(d) }& 2 & $\left\{ \left(\begin{array}{lr}
                                                    1 & \lambda \\
                                                    0 & 1 \\
                                           \end{array} \right), I, (m,n) \right\}, \lambda \neq 0, n \neq 0$ & $Nil^4$  \\ \hline
{\rule [-3mm]{0mm}{8mm}(e) }& 2 &  $\left\{ \left(\begin{array}{lr}
                                                    -1 & \lambda \\
                                                    0 & -1 \\
                                           \end{array} \right), I, (m,n) \right\}, \lambda \neq 0$ & $Nil^3 \times \mathbb{R}$ \\ \hline
{\rule [-3mm]{0mm}{8mm}(f) }& 2 & $\left\{ \left(\begin{array}{lr}
                                                    1 & \lambda \\
                                                    0 & 1 \\
                                           \end{array} \right), -I, (m,n) \right\}, \lambda \neq 0$ & $Nil^3 \times \mathbb{R}$ \\ \hline
{\rule [-3mm]{0mm}{8mm}(g) }& 2 & $\{C, I, (m,n)\}, |tr C| \ge 3, C \in 
SL_2 (\mathbb{Z}) $& $Sol^3 \times \mathbb{R}$ \\ \hline
{\rule [-3mm]{0mm}{8mm}(h) }& 2 & $\{C, -I, (m,n)\}, tr C \ge 3, C \in 
SL_2 (\mathbb{Z}) $& $Sol^3 \times \mathbb{R}$ \\ \hline
\end{tabular}
}

\medskip
\noindent
The given matrices describe the monodromy corresponding to the two generators 
of $\pi_{1}(T^{2})=\bZ^{2}$, and the pairs of integers $(m,n)$ 
represent the Euler class.

\newpage

\bibliographystyle{amsplain}

\begin{thebibliography}{999}

\bibitem{Bande1}
G.~Bande,
{\sl Formes de contact g\'en\'eralis\'e, couples de contact et couples
contacto-symplectiques},
Th\`ese de Doctorat, Universit\'e de Haute Alsace, Mulhouse 2000.

\bibitem{Bande2}
G.~Bande,
{\it Couples contacto-symplectiques},
Trans.~Amer.~Math.~Soc.~{\bf 355} (2003), 1699--1711.

\bibitem{BGK}
G.~Bande, P.~Ghiggini, D.~Kotschick, {\it Stability theorems for 
symplectic and contact pairs}, Int.~Math.~Res.~Not.~(to appear).

\bibitem{BH}
G.~Bande, A.~Hadjar, {\it Contact pairs}, Tohoku Math.~J.~(to appear).

\bibitem{BPV}
W.~Barth, C.~Peters, A.~Van~de~Ven, {\sl Compact Complex Surfaces},
Springer-Verlag, Berlin 1984.

\bibitem{BW}
W.~M.~Boothby, H.~C.~Wang, {\it On contact manifolds},
Ann.~of Math.~{\bf 68} (1958), 721--734.

\bibitem{CG}
G.~Cairns, E.~Ghys, {\it Totally geodesic foliations on 
$4$-manifolds}, J.~Differential Geometry {\bf 23} (1986), 241--254.

\bibitem{Geiges}
H.~Geiges,
{\it Symplectic structures on $T^{2}$-bundles over $T^{2}$}, 
Duke Math.~J.~{\bf 67} (1992), 539--555.

\bibitem{Geiges2}
H.~Geiges,
{\it Symplectic couples on $4$-manifolds}, 
Duke Math.~J.~{\bf 85} (1996), 701--711.

\bibitem{God}
C.~Godbillon, {\sl Feuilletages}, Birkh\"auser Verlag 1991.

\bibitem{Gompf}
R.~E.~Gompf, {\it A new construction of symplectic manifolds},
Ann.~of Math.~{\bf 142} (1995), 527--595.

\bibitem{Hamilton}
M.~Hamilton, {\sl Bi-Lagrangian structures on closed manifolds}, 
Diplomarbeit M\"unchen 2004.

\bibitem{Hil}
J.~A.~Hillman, {\it Flat $4$-manifold groups}, New Zealand 
J.~Math.~{\bf 24} (1995), 29--40.

\bibitem{H}
J.~A.~Hillman, {\sl Four-manifolds, geometries and knots}, Geometry 
$\&$ Topology Monographs, vol.~{\bf 5}, 2002.

\bibitem{Ko}
D.~Kotschick, {\it Orientations and geometrisations of compact complex 
surfaces}, Bull.~London Math.~Soc.~{\bf 29} (1997), 145--149.

\bibitem{K}
D.~Kotschick, {\it On products of harmonic forms}, 
Duke Math.~J.~{\bf 107} (2001), 521--531.

\bibitem{KM}
D.~Kotschick, S.~Morita, {\it Signatures of foliated surface bundles and
the symplectomorphism groups of surfaces}, Topology (to appear).

\bibitem{M}
V.~S.~Matveev, {\it Hyperbolic manifolds are geodesically rigid}, 
Invent.~Math.~{\bf 151} (2003), 579--609. 

\bibitem{Matveev}
V.~S.~Matveev, {\it GF2003 Problem}, Kyoto 2003.

\bibitem{Moser}
J.~Moser,
{\it On the volume elements on a manifold},
Trans.~Amer.~Math.~Soc.~{\bf 120}
(1965), 286--294.

\bibitem{P}
J.~S.~Pasternack, {\it Foliations and compact Lie group actions},
Comment.~Math.~Helv.~{\bf 46} (1971), 467--477.

\bibitem{SF}
K.~Sakamoto, S.~Fukuhara, {\it Classification of $T^{2}$-bundles over 
$T^{2}$}, Tokyo J.~Math.~{\bf 6} (1983), 311--327.

\bibitem{Sh}
I.~H.~Shavel, {\it A class of algebraic surfaces of general type 
constructed from quaternion algebras}, Pacific J.~Math.~{\bf 76} 
(1978), 221-- 245.

\bibitem{SB}
N.~I.~Shepherd-Barron,
{\it Infinite generation for rings of symmetric tensors}, 
Math.~Res.~Letters {\bf 2} (1995), 125--128.

\bibitem{Ue}
M.~Ue, {\it Geometric $4$-manifolds in the sense of Thurston and 
Seifert $4$-manifolds I}, J.~Math.~Soc.~Japan {\bf 42} (1990), 511--540.

\bibitem{wagner}
M.~Wagner, {\sl \"Uber die Klassifikation flacher Riemannscher 
Mannigfaltigkeiten}, Diplomarbeit Basel 1997.

\bibitem{W1}
C.~T.~C.~Wall, {\it Geometries and geometric structures in real dimension 
$4$ and in complex dimension $2$}, in {\sl Geometry and Topology}, 
ed.~J.~Alexander and J.~Harer, Springer Lecture Notes in Math.~{\bf 1167}, Springer 
Verlag 1985.

\bibitem{W2}
C.~T.~C.~Wall, 
{\it Geometric structures on compact complex analytic surfaces}, 
Topology {\bf 25} (1986), 119--153.

\end{thebibliography}


\end{document}